% With \loadreferencesfalse, the file makes its own counters and saves them to "references.tex".
% With \loadreferencestrue, the file loads the counters saved in "references.tex".
%
\newif\ifloadreferences\loadreferencestrue
%
%%%%%%%%%%%%%%%%%%%%%%%%%%%%%%%%%%%%%%%%%%%%%%%%%%%%%%%%%%%%%%%%%%%%%%%%%%%%%%%%%%%%%%%%%%%%%%%%%%%%%%%%%%%%%%%%%%%%%%%
%
% 0: Preliminaries.
%
%%%%%%%%%%%%%%%%%%%%%%%%%%%%%%%%%%%%%%%%%%%%%%%%%%%%%%%%%%%%%%%%%%%%%%%%%%%%%%%%%%%%%%%%%%%%%%%%%%%%%%%%%%%%%%%%%%%%%%%
%
%
\let\myfrac=\frac%
\input eplain %
\let\frac=\myfrac%
\let\myfootnote=\footnote%
\input amstex \input epsf %
\let\footnote=\myfootnote%
%
% Here we load the functions permitting us to use "amsmath" without using "amsppt".
%
\loadeufm\loadmsam\loadmsbm\message{symbol names}\UseAMSsymbols\message{,}%
\magnification 1200 %
\font\myfontdefault=cmr10%
\newif\ifmakebiblio%
\newif\ifinappendices%
\newif\ifundefinedreferences%
\newif\ifchangedreferences%
\makebibliofalse%
\undefinedreferencesfalse%
\changedreferencesfalse%
%
%%%%%%%%%%%%%%%%%%%%%%%%%%%%%%%%%%%%%%%%%%%%%%%%%%%%%%%%%%%%%%%%%%%%%%%%%%%%%%%%%%%%%%%%%%%%%%%%%%%%%%%%%%%%%%%%%%%%%%%
%
% 1: Abstract machinery.
%
% Here we define the macro "makecounter", which gives functionality to the counters defined presently. In order
% to implement it, it is necessary to provisionally change the catcodes.
%
%%%%%%%%%%%%%%%%%%%%%%%%%%%%%%%%%%%%%%%%%%%%%%%%%%%%%%%%%%%%%%%%%%%%%%%%%%%%%%%%%%%%%%%%%%%%%%%%%%%%%%%%%%%%%%%%%%%%%%%
%
\def\setcatcodes{\catcode`\!=0 \catcode`\\=11}%
{\global\let\noe=\noexpand%
\catcode`\@=11 \catcode`\_=11 \setcatcodes%
!global!def!_@@internal@@makeref#1{%
!global!expandafter!def!csname #1ref!endcsname##1{%
!csname _@#1@##1!endcsname%
!expandafter!ifx!csname _@#1@##1!endcsname!relax%
    !write16{#1 ##1 not defined - run saving references}%
    !undefinedreferencestrue%
!fi}}%
!global!def!_@@internal@@makelabel#1{%
!global!expandafter!def!csname #1label!endcsname##1{%
!edef!temptoken{!csname #1info!endcsname}%
!ifloadreferences%
!expandafter!ifx!csname _@#1@##1!endcsname!relax%
!write16{#1 ##1 not hitherto defined - rerun saving references}%
!changedreferencestrue%
!else%
!expandafter!ifx!csname _@#1@##1!endcsname!temptoken%
!else%
!write16{#1 ##1 reference has changed - rerun saving references}%
!changedreferencestrue%
!fi%
!fi%
!else%
!expandafter!edef!csname _@#1@##1!endcsname{!temptoken}%
!edef!textoutput{!write!references{\global\def\_@#1@##1{!temptoken}}}%
!textoutput%
!fi}}%
!global!def!makecounter#1{!_@@internal@@makelabel{#1}!_@@internal@@makeref{#1}}%
!unsetcatcodes%
}
%
%%%%%%%%%%%%%%%%%%%%%%%%%%%%%%%%%%%%%%%%%%%%%%%%%%%%%%%%%%%%%%%%%%%%%%%%%%%%%%%%%%%%%%%%%%%%%%%%%%%%%%%%%%%%%%%%%%%%%%%
%
% 2: Counters.
%
% Here we define the various counters.
%
%%%%%%%%%%%%%%%%%%%%%%%%%%%%%%%%%%%%%%%%%%%%%%%%%%%%%%%%%%%%%%%%%%%%%%%%%%%%%%%%%%%%%%%%%%%%%%%%%%%%%%%%%%%%%%%%%%%%%%%
%
\def\turnintolatin#1{\ifcase #1 _\or i\or ii\or iii\or iv\or v\or vi\or vii\or viii\or ix\or x\or xi\or xii\or xiii\or xiv\or xv\or xvi\or xvii\or xviii\or xix\or xx\or xxi\or xxii\or xxiii\or xxiv\or xxv\or xxvi\fi}%
\def\alphanum#1{\ifcase #1 _\or A\or B\or C\or D\or E\or F\or G\or H\or I\or J\or K\or L\or M\or N\or O\or P\or Q\or R\or S\or T\or U\or V\or W\or X\or Y\or Z\fi}%
\newwrite\references%
\ifloadreferences{\catcode`\@=11 \catcode`\_=11 \global\def\_@citation@AndEtAl{1}
\global\def\_@citation@BarbotBeguinZeghib{2}
\global\def\_@citation@BarbotFillastre{3}
\global\def\_@citation@BonMonSch{4}
\global\def\_@citation@Danciger{5}
\global\def\_@citation@FillSeppi{6}
\global\def\_@citation@FillSmi{7}
\global\def\_@citation@Goldman{8}
\global\def\_@citation@LabourieIII{9}
\global\def\_@citation@Labourie{10}
\global\def\_@citation@LabourieII{11}
\global\def\_@citation@Mess{12}
\global\def\_@citation@Schlenker{13}
\global\def\_@citation@Tamburelli{14}
\global\def\_@citation@Tromba{15}
\global\def\_@subhead@Overview{1}
\global\def\_@subhead@DefinitionsAndStatementOfResult{2}
\global\def\_@eqn@CocycleRelation{\relax \unhbox \voidb@x \hbox {{\relax \tenrm (1)}}}
\global\def\_@eqn@GroupLawOfTwistedProduct{\relax \unhbox \voidb@x \hbox {{\relax \tenrm (2)}}}
\global\def\_@rmk@StructureOfK{\relax \unhbox \voidb@x \hbox {2.1}}
\global\def\_@proc@MainTheorem{2.1}
\global\def\_@rmk@MainTheoremI{\relax \unhbox \voidb@x \hbox {2.2}}
\global\def\_@rmk@MainTheoremII{\relax \unhbox \voidb@x \hbox {2.3}}
\global\def\_@subhead@ProofOfResult{3}
\global\def\_@rmk@TraceFreeTensors{\relax \unhbox \voidb@x \hbox {3.4}}
\global\def\_@eqn@DefinitionOfF{\relax \unhbox \voidb@x \hbox {{\relax \tenrm (3)}}}
\global\def\_@proc@BonMonSch{3.1}
\global\def\_@eqn@FormulaForDerivative{\relax \unhbox \voidb@x \hbox {{\relax \tenrm (4)}}}
\global\def\_@proc@FormulaForDerivative{3.2}
\global\def\_@eqn@OrthogonalDecomposition{\relax \unhbox \voidb@x \hbox {{\relax \tenrm (5)}}}
\global\def\_@proc@OrthogonalDecomposition{3.3}
\global\def\_@eqn@LaplacianEquation{\relax \unhbox \voidb@x \hbox {{\relax \tenrm (6)}}}
\global\def\_@eqn@DefinitionOfPhi{\relax \unhbox \voidb@x \hbox {{\relax \tenrm (7)}}}
\global\def\_@eqn@IntegrabilityCondition{\relax \unhbox \voidb@x \hbox {{\relax \tenrm (8)}}}
\global\def\_@eqn@TheSumIsADerivative{\relax \unhbox \voidb@x \hbox {{\relax \tenrm (9)}}}
\global\def\_@eqn@DefinitionOfX{\relax \unhbox \voidb@x \hbox {{\relax \tenrm (10)}}}
\global\def\_@eqn@FunctionsCoincideAtInfinity{\relax \unhbox \voidb@x \hbox {{\relax \tenrm (11)}}}
\global\def\_@subhead@Bibliography{A}
 }%
\else{\openout\references=references.tex }%
\fi%
%
% A: Headings.
%
\newcount\headno%
\global\headno=0%
\def\headinfo{\ifinappendices\alphanum\headno\else\the\headno\fi}%
\def\nextheadno{\global\advance\headno by 1 \global\subheadno=0 \global\procno=0 \global\eqnno=0 \headinfo}%
\makecounter{head}%
%
% B: Subheadings.
%
\newcount\subheadno%
\global\subheadno=0%
\def\subheadinfo{\ifinappendices\alphanum\subheadno\else\the\subheadno\fi}%
\def\nextsubheadno{\global\advance\subheadno by 1 \global\procno=0 \subheadinfo}%
\makecounter{subhead}%
%
% C: Proclaims (Theorems, Propositions, Lemmas, Corollories, Definitions).
%
\newcount\procno%
\global\procno=0%
\def\procinfo{\subheadinfo.\the\procno}%
\def\nextprocno{\global\advance\procno by 1 \procinfo}%
\makecounter{proc}%
%
% D: Figures.
%
\newcount\figno%
\global\figno=0%
\def\figinfo{\subheadinfo.\the\figno}%
\def\nextfigno{\global\advance\figno by 1 \figinfo}%
\makecounter{fig}%
%
% E: Equations.
%
\newcount\eqnno%
\global\eqnno=0%
\def\eqninfo{\text{{\rm (\the\eqnno)}}}%
\def\nexteqnno[#1]{\global\advance\eqnno by 1 \eqninfo\hbox{\eqnlabel{#1}}}%
\makecounter{eqn}%
%
% F: Remarks.
%
\newcount\rmkno%
\global\rmkno=0%
\def\rmkinfo{\text{\subheadinfo.\the\rmkno}}%
\def\nextrmkno[#1]{\global\advance\rmkno by 1 \rmkinfo\hbox{\rmklabel{#1}}}%
\makecounter{rmk}%
%
%%%%%%%%%%%%%%%%%%%%%%%%%%%%%%%%%%%%%%%%%%%%%%%%%%%%%%%%%%%%%%%%%%%%%%%%%%%%%%%%%%%%%%%%%%%%%%%%%%%%%%%%%%%%%%%%%%%%%%%
%
% 3: Citations.
%
% Citations are treated as a special type of counter.
%
%%%%%%%%%%%%%%%%%%%%%%%%%%%%%%%%%%%%%%%%%%%%%%%%%%%%%%%%%%%%%%%%%%%%%%%%%%%%%%%%%%%%%%%%%%%%%%%%%%%%%%%%%%%%%%%%%%%%%%%
%
\def\gobbleeight#1#2#3#4#5#6#7#8{}%
\newcount\citationno%
\global\citationno=0%
\def\citationinfo{\the\citationno}%
\makecounter{citation}%
\newwrite\biblio%
\def\newref#1#2{%
\def\temptext{#2}%
\edef\bibliotextoutput{\expandafter\gobbleeight\meaning\temptext}%
\global\advance\citationno by 1\citationlabel{#1}%
\ifmakebiblio%
    \edef\fileoutput{\write\biblio{\noindent\hbox to 0pt{\hss$[\the\citationno]$}\hskip 0.2em\bibliotextoutput\medskip}}%
    \fileoutput%
\fi}%
\def\cite#1{%
$[\citationref{#1}]$%
\ifmakebiblio%
    \edef\fileoutput{\write\biblio{#1}}%
    \fileoutput%
\fi%
}%
%
%%%%%%%%%%%%%%%%%%%%%%%%%%%%%%%%%%%%%%%%%%%%%%%%%%%%%%%%%%%%%%%%%%%%%%%%%%%%%%%%%%%%%%%%%%%%%%%%%%%%%%%%%%%%%%%%%%%%%%%
%
% 4: Formatting.
%
%%%%%%%%%%%%%%%%%%%%%%%%%%%%%%%%%%%%%%%%%%%%%%%%%%%%%%%%%%%%%%%%%%%%%%%%%%%%%%%%%%%%%%%%%%%%%%%%%%%%%%%%%%%%%%%%%%%%%%%
%
\let\mypar=\par%
\edef\Pagetitle={Blank}\headline={\hfil\Pagetitle\hfil}%
\edef\Pagefooter={Blank}\footline={\hfil\Pagefooter\hfil}%
%
% A: Move to next odd page.
%
\newcount\showpagenumflag%
\global\showpagenumflag=0 %
\def\nextoddpage%
{\newpage\ifodd\pageno%
\else\global\showpagenumflag=0 %
\null\vfil\eject%
\global\showpagenumflag=1 %
\fi}%
%
% B: Headings.
%
\font\headfont=cmb12%
\def\newhead#1[#2]%
{\ifhmode\mypar\fi%
\ifnum\headno=0 \else\goodbreak\bigskip\fi%
{\headfont\noindent\nextheadno\ - #1.}\headlabel{#2}%
\nobreak\medskip}%
%
% C: Subheadings.
%
\def\newsubhead#1[#2]%
{\ifhmode\mypar\fi%
\ifnum\subheadno=0 \else\goodbreak\medskip\fi%
{\bf\noindent\nextsubheadno\ - #1.\ }\subheadlabel{#2}}%
%
% D: Proclaims.
%
\newif\ifinproclaim%
\global\inproclaimfalse%
\def\proclaim#1{%
\goodbreak\medskip
\bgroup\inproclaimtrue%
\noindent{\bf #1}%
\nobreak\medskip\sl}%
\def\noskipproclaim#1{%
\goodbreak\medskip%
\bgroup\inproclaimtrue%
\noindent{\bf #1}\nobreak\sl}%
\def\endproclaim{\mypar\egroup\nobreak\medskip\ignorespaces}%
%
% E: Figures.
%
% The macro "\makelabelgrid" is a useful utility for guiding the positioning of labels is figures.
%
\newcount\xpos\newcount\ypos
\def\makelabelgrid{%
\xpos=-5 \ypos=-5 %
\loop\ifnum\xpos<6 %
{\loop\ifnum\ypos<6 %
\def\labeltext{x}%
\ifnum\xpos=0\def\labeltext{+}\fi%
\ifnum\ypos=0\def\labeltext{+}\fi%
\placelabel[\xpos][\ypos]{\labeltext}%
\advance\ypos by 1 %
\repeat}%
\advance\xpos by 1 %
\repeat}%
\def\placelabel[#1][#2]#3{{%
\setbox10=\hbox{\raise #2cm \hbox{\hskip #1cm #3}}%
\ht10=0pt \dp10=0pt \wd10=0pt \box10}}%
%
%
%
% F: Items.
%
\def\myitem#1{\noindent\hbox to .5cm{\hfill#1\hss}}%
%
% G: Right justification.
%
%
%
%%%%%%%%%%%%%%%%%%%%%%%%%%%%%%%%%%%%%%%%%%%%%%%%%%%%%%%%%%%%%%%%%%%%%%%%%%%%%%%%%%%%%%%%%%%%%%%%%%%%%%%%%%%%%%%%%%%%%%%
%
% 5: Special fonts.
%
%%%%%%%%%%%%%%%%%%%%%%%%%%%%%%%%%%%%%%%%%%%%%%%%%%%%%%%%%%%%%%%%%%%%%%%%%%%%%%%%%%%%%%%%%%%%%%%%%%%%%%%%%%%%%%%%%%%%%%%
%
% A: The "mathsf" font is not defined in Plain.
%
%
\font\sansseriften=cmss10%
\font\sansserifseven=cmss7%
\font\sansseriffive=cmss5%
\newfam\sansseriffam%
\textfont\sansseriffam=\sansseriften%
\scriptfont\sansseriffam=\sansserifseven%
\scriptscriptfont\sansseriffam=\sansseriffive%
\def\mathsf{\fam\sansseriffam}%
%
%
% B: The "mathbf" font is not defined in Plain.
%
\font\boldten=cmb10%
\font\boldseven=cmb7%
\font\boldfive=cmb5%
\newfam\mathboldfam%
\textfont\mathboldfam=\boldten%
\scriptfont\mathboldfam=\boldseven%
\scriptscriptfont\mathboldfam=\boldfive%
\def\mathbf{\fam\mathboldfam}%
%
%
% C: Here we define the macros "mathi" and "mathj". This is not really necessary, since the macros "imath" and
% "jmath" perform the same function. Still, it makes for an interesting exercise.
%
\font\mycmmiten=cmmi10%
\font\mycmmiseven=cmmi7%
\font\mycmmifive=cmmi5%
\newfam\mycmmifam%
\textfont\mycmmifam=\mycmmiten%
\scriptfont\mycmmifam=\mycmmiseven%
\scriptscriptfont\mycmmifam=\mycmmifive%
\def\hexa#1{\ifcase #1 0\or 1\or 2\or 3\or 4\or 5\or 6\or 7\or 8\or 9\or A\or B\or C\or D\or E\or F\fi}%
\mathchardef\mathi="7\hexa\mycmmifam7B%
\mathchardef\mathj="7\hexa\mycmmifam7C%
%
% D: Here we define a few Hebrew letters: "mybeth", "mygimmel" and "mydaleth".
%
\font\mymsbmten=msbm10 at 8pt%
\font\mymsbmseven=msbm7 at 5.6pt%6
\font\mymsbmfive=msbm5 at 4pt%
\newfam\mymsbmfam%
\textfont\mymsbmfam=\mymsbmten%
\scriptfont\mymsbmfam=\mymsbmseven%
\scriptscriptfont\mymsbmfam=\mymsbmfive%
\mathchardef\mybeth="7\hexa\mymsbmfam69%
\mathchardef\mygimmel="7\hexa\mymsbmfam6A%
\mathchardef\mydaleth="7\hexa\mymsbmfam6B%
%
%%%%%%%%%%%%%%%%%%%%%%%%%%%%%%%%%%%%%%%%%%%%%%%%%%%%%%%%%%%%%%%%%%%%%%%%%%%%%%%%%%%%%%%%%%%%%%%%%%%%%%%%%%%%%%%%%%%%%%%
%
% 6: Mathematics operators and symbols.
%
%%%%%%%%%%%%%%%%%%%%%%%%%%%%%%%%%%%%%%%%%%%%%%%%%%%%%%%%%%%%%%%%%%%%%%%%%%%%%%%%%%%%%%%%%%%%%%%%%%%%%%%%%%%%%%%%%%%%%%%
%
\def\proof{{\noindent\bf Proof:\ }}%
\def\remark[#1]{{\noindent\bf Remark \nextrmkno[#1].}}%
\def\qed{~$\square$}%
\def\makeop#1{\global\expandafter\def\csname op#1\endcsname{{\text{#1}}}}%
\def\makeopsmall#1{\global\expandafter\def\csname op#1\endcsname{{\text{\lowercase{#1}}}}}%
%
% A: Set Theory.
%
\def\munion{\mathop{\cup}}%
%
%
% B: Point set topology.
%
\makeop{Ext}%
\makeop{Int}%
\makeop{Dist}%
\makeop{Diam}%
\makeop{Length}%
%
% C: Sequences.
%
%
%
%
\def\mlim{\mathop{{\text{Lim}}}}%
%
%
%
%
%
% D: Linear Algebra.
%
\makeop{Dim}%
\makeop{Ker}%
\makeop{Coker}%
\makeop{Tr}%
\makeop{Adj}%
\makeop{Det}%
\makeop{End}%
\makeop{Lin}%
\makeop{Symm}%
\makeop{Mult}%
%
% E: Basic calculus.
%
\makeop{dx}%
\makeop{dy}%
\makeop{dz}%
\makeop{dt}%
\makeop{dVol}%
\makeop{dArea}%
\makeop{Supp}%
\makeop{Hess}%
\makeop{Lip}%
%
% F: Complex Numbers.
%
\makeop{Re}%
\makeop{Im}%
\makeop{Arg}%
\makeop{Log}%
\makeop{Exp}%
%
% G: Trigonometry.
%
\makeopsmall{Cos}%
\makeopsmall{Sin}%
\makeopsmall{Tan}%
\makeopsmall{Sec}%
\makeopsmall{Cosec}%
\makeopsmall{Cot}%
\makeopsmall{ArcCos}%
\makeopsmall{ArcSin}%
\makeopsmall{ArcTan}%
\makeopsmall{ArcSec}%
\makeopsmall{ArcCosec}%
\makeopsmall{ArcCot}%
%
% H: Hyperbolic Trigonometry.
%
\makeopsmall{Cosh}%
\makeopsmall{Sinh}%
\makeopsmall{Tanh}%
\makeopsmall{ArcCosh}%
\makeopsmall{ArcSinh}%
\makeopsmall{ArcTanh}%
%
% I: Differential and Riemannian Geometry.
%
\makeop{Vol}%
\makeop{Area}%
\makeop{Riem}%
\makeop{Ric}%
\makeop{Scal}%
\makeop{Euc}%
\makeop{Imm}%
\makeop{Emb}%
%
% J: Lie Groups.
%
\makeop{Id}%
\makeop{Ad}%
\makeop{O}%
\makeop{SO}%
\makeop{SL}%
\makeop{GL}%
\makeop{Conf}%
\makeop{Homeo}%
\makeop{Diff}%
\makeop{Isom}%
%
% K: Functional Analysis.
%
\makeop{Ind}%
\makeop{Sig}%
\makeop{Spec}%
%
% L: Other.
%
\makeop{Conv}%
\makeop{Max}%
\makeop{Min}%
\makeop{Mod}%
\makeop{Deg}%
\makeop{loc}%
%
%%%%%%%%%%%%%%%%%%%%%%%%%%%%%%%%%%%%%%%%%%%%%%%%%%%%%%%%%%%%%%%%%%%%%%%%%%%%%%%%%%%%%%%%%%%%%%%%%%%%%%%%%%%%%%%%%%%%%%%
%
% 7: Redundant Material.
%
%%%%%%%%%%%%%%%%%%%%%%%%%%%%%%%%%%%%%%%%%%%%%%%%%%%%%%%%%%%%%%%%%%%%%%%%%%%%%%%%%%%%%%%%%%%%%%%%%%%%%%%%%%%%%%%%%%%%%%%
%
% This file contains redundant functionality for constructing a bibliography. Although it is not used, it may one day
% prove useful, and so I leave it here.
%
% A: Before the citations, the file should contain:
%
% \newif\ifmakebiblio
%
% followed by either \makebibliotrue or \makebibliofalse.
%
% (note that, for conveniance, the commands \newif\ifmakebiblio and \makebibliofalse have been included at the
% beginning of this preamble. These commands should be removed before making the changes outlined here.
%
% B: Immediately before the first citation, the file should contain the following instructions:
%
% \ifmakebiblio%
% \openout\biblio=biblio.tex %
% \edef\fileoutput{\write\biblio{\bgroup\leftskip=2em}}%
% \fileoutput%
% \fi%
%
% C: Immediately after the last citation, the file should contain the following instruction:
%
% \ifmakebiblio%
% {\edef\fileoutput{\write\biblio{\egroup}}%
% \fileoutput}%
% \fi%
%
 %
%
% In order to use "dvips", enter:
% dvips -N0 -Z0 -K0 [whatever.dvi] -o [whatever.ps]
%
%%%%%%%%%%%%%%%%%%%%%%%%%%%%%%%%%%%%%%%%%%%%%%%%%%%%%%%%%%%%%%%%%%%%%%%%%%%%%%%%%%%%%%%%%%%%%%%%%%%%%%%%%%%%%%%%%%%%%%%
%
% 1: The Paper.
%
%%%%%%%%%%%%%%%%%%%%%%%%%%%%%%%%%%%%%%%%%%%%%%%%%%%%%%%%%%%%%%%%%%%%%%%%%%%%%%%%%%%%%%%%%%%%%%%%%%%%%%%%%%%%%%%%%%%%%%%
%
\def\Pagetitle{\hfil}
\def\Pagefooter{\hfil}

\makeop{dS}
\makeop{ad}
\makeop{AdS}
\makeop{PSO}
\makeop{I}
\makeop{II}
\makeop{III}
\makeop{Width}
\makeop{hyp}
\makeop{rep}
\makeop{T}
\makeop{SGr}
\makeop{Coth}
\makeop{GHMC}
\makeop{EC}
\makeop{Lam}
\makeop{PSL}
\makeop{T}
\def\Thyp{{\text{T}_\ophyp}}
\def\Trep{{\text{T}_\oprep}}

\font\tablefont=cmr7
\newif\ifshowaddress\showaddresstrue
\null \vfill
\def\centre{\rightskip=0pt plus 1fil \leftskip=0pt plus 1fil \spaceskip=.3333em \xspaceskip=.5em \parfillskip=0em \parindent=0em}%
\def\textmonth#1{\ifcase#1\or January\or Febuary\or March\or April\or May\or June\or July\or August\or September\or October\or November\or December\fi}
\font\abstracttitlefont=cmr10 at 14pt {\abstracttitlefont\centre A note on invariant constant curvature immersions in Minkowski space.\par}
\bigskip
{\centre 15th August 2018\par}
%{\centre \the\day\ \textmonth\month\ \the\year\par}
\bigskip
{\centre Fran\c{c}ois Fillastre\footnote{${}^1$}{{\tablefont Universit\'e de Cergy-Pontoise, UMR CNRS 8088, Cergy-Pontoise 95000, France\hfill}},
Graham Smith\footnote{${}^2$}{{\tablefont Instituto de Matem\'atica, UFRJ, Av. Athos da Silveira Ramos 149, Centro de Tecnologia - Bloco C, Cidade Universit\'aria - Ilha do Fund\~ao, Caixa Postal 68530, 21941-909, Rio de Janeiro, RJ - BRAZIL\hfill}}\par}
\bigskip
\noindent{\bf Abstract:~}Let $S$ be a compact, orientable surface of hyperbolic type. Let $(k_+,k_-)$ be a pair of negative numbers and let $(g_+, g_-)$ be a pair of marked metrics over $S$ of constant curvature equal to $k_+$ and $k_-$ respectively. Using a functional introduced by Bonsante, Mondello \& Schlenker, we show that there exists a unique affine deformation $\Gamma:=(\rho,\tau)$ of a Fuchsian group such that $(S,g_+)$ and $(S, g_-)$ embed isometrically as locally strictly convex Cauchy surfaces in the future and past complete components respectively of the quotient by $\Gamma$ of an open subset $\Omega$ of Minkowski space.
\par
Such quotients are known as Globally Hyperbolic, Maximal, Cauchy compact Min\-kow\-ski spacetimes and are naturally dual to the half-pipe spaces introduced by Danciger. When translated into this latter framework, our result states that there exists a unique, marked, quasi-Fuchsian half-pipe space in which $(S, g_+)$ and $(S, g_-)$ are realised as the third fundamental forms of future- and past-oriented, locally strictly convex graphs.
\bigskip
\noindent{\bf Classification AMS~:~}30F60, 53C50
%
% 30F60, Teichmueller Theory
% 53C50, Lorentz manifold
%
\par
\vfill
\eject
%\nextoddpage
%
\global\pageno=1
\myfontdefault
\def\Pagetitle{\hfil A note on invariant constant curvature immersions in Minkowski space.\hfil}
\def\Pagefooter{\hfil{\myfontdefault\folio}\hfil}
\catcode`\@=11
\def\triplealign#1{\null\,\vcenter{\openup1\jot \m@th %
\ialign{\strut\hfil$\displaystyle{##}\quad$&$\displaystyle{{}##}$\hfil&$\displaystyle{{}##}$\hfil\crcr#1\crcr}}\,}
\def\multiline#1{\null\,\vcenter{\openup1\jot \m@th %
\ialign{\strut$\displaystyle{##}$\hfil&$\displaystyle{{}##}$\hfil\crcr#1\crcr}}\,}
\catcode`\@=12
\newref{AndEtAl}{Andersson L., Barbot T., Benedetti R., Bonsante F., Goldman W. M., Labourie F., Scannell K. P., Schlenker J. M., Notes on: ``Lorentz spacetimes of constant curvature'', {\sl Geom. Dedicata}, {\bf 126}, (2007), 47--70}
\newref{BarbotBeguinZeghib}{Barbot T., B\'eguin F., Zeghib A., Prescribing Gauss curvature of surfaces in 3-dimensional spacetimes: application to the Minkowski problem in the Minkowski space, {\sl Ann. Inst. Fourier}, {\bf 61}, (2011), no. 2, 511--591}
\newref{BarbotFillastre}{Barbot T., Fillastre F., Quasi-Fuchsian co-Minkowski manifolds, arXiv:1801.10429}
\newref{BonMonSch}{Bonsante F., Mondello G., Schlenker J. M., A cyclic extension of the earthquake flow {II}, {\sl Ann. Sci. \'Ec. Norm. Sup\'er.}, {\bf 48}, (2015), no. 4, 811--859}
\newref{Danciger}{Danciger J., A Geometric transition from hyperbolic to anti de Sitter geometry, {\sl Geom. Topol.}, {\bf 17}, (2013), no. 5, 3077--3134}
\newref{FillSeppi}{Fillastre F., Seppi A., Spherical, hyperbolic and other projective geometries, in {\sl Eighteen essays on non-commutative geometry}, IRMA, {\sl Math. Theor. Phys.}, {\bf 29}, Eur. Math. Soc., Z\"urich, (2019)}
\newref{FillSmi}{Fillastre F., Smith G., Group actions and scattering problems in Teichm\"uller theory, arXiv:1605.04563}
\newref{Goldman}{Goldman W., The symplectic nature of the fundamental group of surfaces, {\sl Adv. Math.}, {\bf 54}, (1984), no. 2, 200--225}
\newref{LabourieIII}{Labourie F., Probl\`eme de Minkowski et surfaces \`a courbure constante dans les vari\'et\'es hyperboliques, {\sl Bull. Soc. math. France}, {\bf 119}, 1991, 307--325}
\newref{Labourie}{Labourie F., Surfaces convexes dans l'espace hyperbolique et $\Bbb{CP}^1$-structures, {\sl J. London Math. Soc.}, {\bf 45}, (1992), no. 3, 549--565}
\newref{LabourieII}{Labourie F., M\'etriques prescrites sur le bord des vari\'et\'es hyperboliques de dimension $3$, {\sl J. Diff. Geom.}, {\bf 35}, (1992), no. 3, 609--626}
\newref{Mess}{Mess G., Lorentz spacetimes of constant curvature, {\sl Geom. Dedicata}, {\bf 126}, (2007), 3--45}
\newref{Schlenker}{Schlenker J. M., Hyperbolic manifolds with convex boundary, {\sl Inventiones math.}, {\bf 163}, (2006), 109--169}
\newref{Tamburelli}{Tamburelli A., Prescribing metrics on the boundary of anti-de Sitter $3$-manifolds, {\sl Int. Math. Res. Not. IMRN}, (2018), no. 5, 1281--1313}
\newref{Tromba}{Tromba A. J., {\sl Teichm\"uller theory in Riemannian geometry}, Lectures in Mathematics ETH Z\"urich, Birkh\"auser Verlag, Basel, (1992)}
\newsubhead{Overview}[Overview]
Let $S$ be a compact, orientable surface of hyperbolic type, and let $(g_+,g_-)$ be a pair of marked metrics over $S$. In \cite{LabourieIII}, Labourie proved that when these metrics have constant curvature equal to $k_\pm\in]-1,0[$, there exists a marked quasi-Fuchsian hyperbolic $3$-manifold $M$, unique up to orientation, into which $(S,g_+)$ and $(S,g_-)$ embed isometrically as locally strictly convex surfaces in each of the components of the complement of its Nielsen kernel. As a consequence, Labourie obtained a $2$-dimensional, real-analytic family of real-analytic parametrisations of the space of marked quasi-Fuchsian hyperbolic $3$-manifolds by two copies of the Teichm\"uller space of $S$. In particular, this family interpolates between the double parametrisation of Ahlfors-Bers and the double parametrisation, conjectured by Thurston, in terms of the intrinsic hyperbolic metrics of the two sides of the convex hull (c.f. Theorem $5.3$ of \cite{FillSmi}).
\par
In \cite{LabourieII}, Labourie extended the existence part of this result to the case where $g_\pm$ have non-constant curvature in $]-1,0[$ and uniqueness was subsequently proven by Schlenker in \cite{Schlenker}. Schlenker proved furthermore that, under similar, though slightly stronger, hypotheses, there also exists a marked quasi-Fuchsian hyperbolic manifold, which is unique up to orientation, in which $(S,g_+)$ and $(S,g_-)$ are realised as the {\sl third} fundamental forms of locally strictly convex immersions in each of the components of the complement of the Nielsen kernel. That is, the initial result of Labourie, which concerns the prescription of first fundamental forms, has a precise analogue concerning the prescription of third fundamental forms. A classical duality argument (c.f. \cite{FillSeppi}) then shows that this latter result is equivalent to the existence of a unique, time-oriented, marked quasi-Fuchsian {\sl de Sitter spacetime} into which $(S,g_+)$ and $(S,g_-)$ embed isometrically as locally strictly convex Cauchy surfaces in its future and past components respectively.
\par
This result was adapted to the anti de Sitter case by Bonsante, Mondello \& Schlenker in \cite{BonMonSch}. To be precise, they show that when $g_+$ and $g_-$ have constant curvature in $]-\infty,-1[$, there exists a marked, time-oriented, quasi-Fuchsian anti de Sitter spacetime into which $(S,g_+)$ and $(S,g_-)$ embed as locally strictly convex Cauchy surfaces in its future and past components respectively. In this case, the same duality argument mentioned above yields the analogous result concerning prescription of third fundamental forms of surfaces in anti de Sitter spacetimes. The existence part of Bonsante, Mondello \& Schlenker's result is extended to the case of metrics of non-constant curvature by Tamburelli in \cite{Tamburelli}. However, even in the case of constant curvature, Bonsante, Mondello \& Schlenker only proved uniqueness when the curvatures $k_\pm$ of $g_\pm$ are related by
$$
(k_+ + 1)(k_- + 1) = 1,
$$
and the problem of uniqueness in its full generality remains unsolved.
\par
In this paper, we will be concerned with the intermediate case between quasi-Fuchsian hyperbolic manifolds and quasi-Fuchsian anti de Sitter manifolds. This is known to be that of quasi-Fuchsian co-Minkowski (or half-pipe) manifolds, introduced by Danciger in \cite{Danciger} (c.f. also the review \cite{BarbotFillastre} of Barbot and the first author). Here, the ambient space carries a natural {\sl degenerate} metric preserved by a natural covariant derivative of constant curvature equal to $-1$. From these properties it follows that any smooth graph is intrinsically hyperbolic. In the present paper, we adapt the above results to the half-pipe case and show that, given two marked metrics $(g_+,g_-)$ on $S$ of constant negative curvature, there exists a unique, marked half-pipe space into which $(S,g_+)$ and $(S,g_-)$ are realised respectively as the third fundamental forms of future- and past-oriented locally strictly convex graphs.
\newsubhead{Definitions and statement of result}[DefinitionsAndStatementOfResult]
By duality, our result has an equivalent statement in terms of Minkowski spacetimes, which arise in the theory of general relativity. This is the framework that we will use throughout the sequel (we refer the reader to \cite{FillSmi} for an in-depth review). A {\sl Minkowski spacetime} is a semi-riemannian manifold which is everywhere locally isometric to $\Bbb{R}^{2,1}$. A smoothly embedded curve in a Minkowski spacetime is said to be {\sl causal} whenever its derivative has non-positive norm-squared at every point. A Minkowski spacetime is itself said to be {\sl causal} whenever it contains no non-trivial, closed, causal curve. A causal spacetime is said to be {\sl globally hyperbolic} whenever it contains a {\sl Cauchy hypersurface}, that is, a smoothly embedded hypersurface that meets every inextensible, causal curve exactly once. A globally hyperbolic spacetime $X$ is said to be {\sl maximal} whenever it cannot be isometrically embedded into a strictly larger globally hyperbolic spacetime $X'$ in such a manner that the Cauchy hypersurfaces of $X$ are mapped to Cauchy hypersurfaces of $X'$. Finally, a globally hyperbolic spacetime is said to be {\sl Cauchy compact} whenever its Cauchy hypersurface, which is unique up to diffeomorphism, is compact. A Minkowski spacetime which possesses all the above properties is said to be {\sl GHMC} (Globally, Hyperbolic, Maximal and Cauchy Compact).
\par
Let $\opGHMC_0[S]$ denote the space of all marked, GHMC Minkowski spacetimes with Cauchy surface diffeomorphic to $S$ where we recall that $S$ is a compact, orientable surface of hyperbolic type. We recall the explicit parametrisation of $\opGHMC_0[S]$ in terms of algebraic data constructed by Mess in \cite{Mess} (c.f. also \cite{AndEtAl} and \cite{FillSmi}). To this end, we first introduce affine deformations. Thus, let $\pi_1(S)$ denote the fundamental group of $S$. A homomorphism $\rho:\pi_1(S)\rightarrow\opSO(2,1)$ is said to be {\sl Fuchsian} whenever its image acts properly discontinuously without torsion on $\Bbb{H}^2$. An {\sl affine deformation} of $S$ is defined to be a pair $(\rho,\tau)$, where $\rho$ is a Fuchsian homomorphism and $\tau:\pi_1(S)\rightarrow\Bbb{R}^{2,1}$ satisfies the {\sl $\rho$-cocycle condition}%
$$
\tau(gh) = \tau(g) + \rho(g)\tau(h),\eqnum{\nexteqnno[CocycleRelation]}
$$
for all $g,h\in\pi_1(S)$. Every affine deformation $(\rho,\tau)$ identifies with a homomorphism $g\mapsto(\rho(g),\tau(g))$ of $\pi_1(S)$ into the twisted product $\opSO(2,1)\ltimes\Bbb{R}^{2,1}$, where the group law of the latter is given by
$$
(g,X)\cdot(h,Y) := (gh, X + gY).\eqnum{\nexteqnno[GroupLawOfTwistedProduct]}
$$
In the sequel, $\rho$ will be called the {\sl Fuchsian component} of the affine deformation, $\tau$ will be called its {\sl cocycle} and the image of the above homomorphism will be denoted by $\rho\ltimes\tau$. Finally, recall that the space $\Trep[S]$ of $\opSO(2,1)$-conjugacy classes of Fuchsian homomorphisms coincides with the Teichm\"uller space of the surface $S$. In like manner, the space of $\opSO(2,1)\ltimes\Bbb{R}^{2,1}$-conjugacy classes of affine deformations identifies with the tangent bundle $\opT\Trep[S]$ of $\Trep[S]$ (c.f. \cite{Goldman}).
\par
Observe that $\opSO(2,1)\ltimes\Bbb{R}^{2,1}$ is the identity component of the group of rigid, affine motions of $\Bbb{R}^{2,1}$. In \cite{Mess} (c.f. also \cite{AndEtAl}), Mess shows that, for every affine deformation $(\rho,\tau)$, there exists a unique pair $K_+$ and $K_-$ of convex subsets of $\Bbb{R}^{2,1}$ such that
\medskip
\myitem{(1)} all supporting normals to $K_\pm$ are either timelike or null;
\medskip
\myitem{(2)} $K_+$ and $K_-$ are respectively future- and past-complete;
\medskip
\myitem{(3)} $\rho\ltimes\tau$ acts properly discontinuously over the interiors of $K_+$ and $K_-$; and
\medskip
\myitem{(4)} $(K_+,K_-)$ is maximal among all such pairs in the sense that if $(K_+',K_-')$ is another pair of convex subsets of $\Bbb{R}^{2,1}$ possessing the above three properties, then $K_+'\subseteq K_+$ and $K_-'\subseteq K_-$.
\medskip
\remark[StructureOfK] It follows from the definition that $K_+$ (resp. $K_-$) is the region in $\Bbb{R}^{2,1}$ bounded below (resp. above) by a strictly convex (resp. concave), Lipschitz graph over $\Bbb{R}^2$.
\medskip
\noindent The interiors of $K_+$ and $K_-$ are referred to respectively as the {\sl Mess' future and past domains} of the affine deformation $(\rho,\tau)$. The quotients $K_\pm/(\rho\ltimes\tau)$ define respectively future-complete and past-complete, marked, GHMC Minkowski spacetimes. Furthermore, two affine deformations yield the same marked, GHMC Minkowski spacetime if and only if they lie in the same conjugacy class, and all future- and past-complete, marked, GHMC Minkowski spacetimes arise in this manner. This concludes the construction of Mess' parametrisation of $\opGHMC_0[S]$ by $\opT\Trep[S]$. In particular, this identification furnishes $\opGHMC_0[S]$ with a real-analytic structure. Furthermore, we see that every marked, GHMC Minkowski spacetime $X$ in fact consists of one future-complete component $X_+$ and one past-complete component $X_-$.
\par
Let $\Thyp[S]$ denote the Teichm\"uller space of marked, hyperbolic metrics over $S$. Using surfaces of constant extrinsic curvature, we define a real-analytic function
$$
\Phi:]-\infty,0[\times]-\infty,0[\times\opGHMC_0[S]\rightarrow\Thyp[S]\times\Thyp[S]
$$
as follows. First, in \cite{BarbotBeguinZeghib}, Barbot, B\'eguin \& Zeghib show that for every GHMC Minkowski spacetime $X$, and for all $k_+,k_-\in]-\infty,0[$, there exists a unique pair $(\Sigma_+,\Sigma_-)$ where
\medskip
\myitem{(1)} $\Sigma_\pm$ is a smooth Cauchy hypersurface in $X_\pm$;
\medskip
\myitem{(2)} $\Sigma_\pm$ is locally strictly convex in the sense that its shape operator is everywhere positive definite; and
\medskip
\myitem{(3)} $\Sigma_\pm$ has constant extrinsic curvature equal to $k_\pm$.
\medskip
\noindent Up to rescaling, the intrinsic metrics of $\Sigma_\pm$ are hyperbolic and the marking on $X$ induces markings on $\Sigma_\pm$. In this manner $\Sigma_\pm$ define points in $\Thyp[S]$, and the function $\Phi$ is given by
$$
\Phi(k_+,k_-,X):=(\Sigma_+,\Sigma_-).
$$
We prove
\proclaim{Theorem \nextprocno}
\noindent For all $(k_+,k_-)\in]-\infty,0[^2$, the map $\Phi(k_+,k_-,\cdot)$ defines a real-analytic diffeomorphism from $\opGHMC_0[S]$ into $\Thyp[S]\times\Thyp[S]$.
\endproclaim
\proclabel{MainTheorem}
\remark[MainTheoremI] In particular, for any pair $(g_+,g_-)$ of marked metrics over $S$ of constant negative curvature, there exists a unique, marked, GHMC Minkowski spacetime $X$ such that $(S,g_+)$ and $(S,g_-)$ embed isometrically as locally strictly convex Cauchy surfaces in its future- and past-complete components respectively.
\medskip
\remark[MainTheoremII] Theorem \procref{MainTheorem}, together with a formal adaptation of the arguments of \cite{LabourieII} and \cite{Tamburelli}, yields {\sl existence} also in the case of variable curvature. That is, for any pair $(g_+,g_-)$ of marked, smooth, negative curvature metrics over $S$, there exists a marked, GHMC Minkowski spacetime $X:=X_+\munion X_-$ such that $(S,g_\pm)$ embeds isometrically into $X_\pm$.
\newsubhead{Proof of result}[ProofOfResult]
Given a hyperbolic metric $h$ over $S$, a {\sl Codazzi field} of $h$ is defined to be a smooth section $M$ of $\opEnd(TS)$ such that
\medskip
\myitem{(1)} $M$ is symmetric with respect to $h$; and
\medskip
\myitem{(2)} $d^\nabla M=0$, where the covariant derivative is also taken with respect to $h$.
\medskip
\remark[TraceFreeTensors] For any section $M$ of $\opEnd(TS)$,
$$
d^\nabla M = (\nabla\cdot MJ)\opdArea_h,
$$
where $\nabla\cdot$ here denotes the divergence operator of $h$, $J$ its complex structure, and $\opdArea_h$ its area form. Furthermore, a section $M$ of $\opEnd(TS)$ is symmetric and trace-free if and only if the section $MJ$ is also symmetric and trace free. In \cite{Tromba}, Tromba studies symmetric, trace-free sections $M$ of $\opEnd(TS)$ which satisfy $\nabla\cdot M=0$. It follows from the preceeding observations that Tromba's formalism transforms into our own by multiplication on the right by $J$.
\medskip
\noindent Recall now from \cite{Labourie} that, given two marked hyperbolic metrics $h$ and $h'$, there exists a unique Codazzi field $B:=B(h',h)$ of $h$ such that
\medskip
\myitem{(1)} $\opDet(B)=1$; and
\medskip
\myitem{(2)} $h(B\cdot,B\cdot)$ coincides with $h'$ as a point in $\Thyp[S]$.
\medskip
\noindent This field is referred to as the {\sl Labourie field} of $h'$ with respect to $h$. Consider now the function $F:\Thyp[S]\times\Thyp[S]\rightarrow\Bbb{R}$ given by
$$
F_{h'}(h) := F(h',h) := \int_\Sigma\opTr(B(h',h))\opdArea_h.\eqnum{\nexteqnno[DefinitionOfF]}
$$
In \cite{BonMonSch}, Bonsante, Mondello \& Schlenker establish the following properties of $F$.
\proclaim{Lemma \nextprocno, Bonsante, Mondello \& Schlenker {\bf \cite{BonMonSch}}}
\myitem{(1)} $F$ is symmetric in $h$ and $h'$.
\medskip
\noindent Furthermore, for all $h'$,
\medskip
\myitem{(2)} $F_{h'}$ is proper;
\medskip
\myitem{(3)} $F_{h'}$ is strictly convex with respect to the Weyl-Peterson metric; and
\medskip
\myitem{(4)} $F_{h'}$ attains a unique minimum at $h=h'$.
\endproclaim
\proclabel{BonMonSch}
\noindent Bonsante, Mondello \& Schlenker also determine a formula for the derivative of this functional. For the reader's convenience, we provide a simpler proof of this relation. First recall that infinitesimal variations of a given hyperbolic metric are given by trace-free Codazzi fields of that metric (c.f. \cite{Tromba}). We now have
\proclaim{Lemma \nextprocno}
\noindent For every trace-free Codazzi field $A$ of $h$,
$$
DF_{h'}(h)A = -\frac{1}{2}\int_\Sigma\opTr(AB(h',h))\opdArea_{h}.\eqnum{\nexteqnno[FormulaForDerivative]}
$$
\endproclaim
\proclabel{FormulaForDerivative}
\proof Denote $B:=B(h',h)$ and identify $h'$ with $h(B\cdot,B\cdot)$. Consider an infinitesimal perturbation of $h$ given by the trace-free Codazzi field $A$, that is
$$
\delta h = h(A\cdot,\cdot).
$$
By definition,
$$\eqalign{
\opTr(A) &= 0,\ \text{and}\cr
d^{\nabla} A&=0,\cr}
$$
where $\nabla$ here denotes the covariant derivative of $h$. Let $\delta B$ be the resulting infinitesimal variation of the Labourie field. We have
$$
h(AB\cdot,B\cdot) + h(\delta B\cdot,B\cdot) + h(B\cdot,\delta B\cdot) = (\Cal{L}_Xh')(\cdot,\cdot),
$$
for some vector field $X$. From this we deduce that
$$
BAB + B\delta B + \delta B B = BM + M^tB,
$$
where $M:=\nabla(BX)$. Composing with $B^{-1}$ and taking the trace then yields
$$
\opTr(\delta B) = -\frac{1}{2}\opTr(AB) + \nabla\cdot(BX).
$$
Finally, since $A$ is trace-free, the resulting infinitesimal variation of the area form $\opdArea_h$ vanishes. It follows that
$$
DF_{h'}(h)A = \int_\Sigma\opTr(\delta B)\opdArea_{h},
$$
and the result now follows by Stokes' Theorem.\qed
\medskip
The following orthogonal decomposition will also prove useful.
\proclaim{Lemma \nextprocno}
\noindent Let $h$ be a hyperbolic metric over $S$. For every Codazzi field $M$ of $h$, there exists a unique smooth function $f:S\rightarrow\Bbb{R}$ and a unique trace-free Codazzi field $A$ of $h$ such that
$$
M = A + \big(f\opId - \opHess_h(f)\big),\eqnum{\nexteqnno[OrthogonalDecomposition]}
$$
where $\opHess_h$ denotes the Hessian with respect to $h$. Furthermore, this decomposition is orthogonal with respect to the $L^2$-norm of $h$.
\endproclaim
\proclabel{OrthogonalDecomposition}
\proof Observe that when \eqnref{OrthogonalDecomposition} is satisfied, we have
$$
(\Delta^h - 2)f = -\opTr(M),\eqnum{\nexteqnno[LaplacianEquation]}
$$
where $\Delta^h$ here denotes the Laplace operator of $h$. However, by the maximum principle, the operator $(\Delta^h-2)$ is injective and, by the Fredholm alternative, it is a linear isomorphism. Uniqueness of $f$, and therefore also of $A$, thus follows. Conversely, if $f$ is the unique solution of \eqnref{LaplacianEquation}, then
$$
A := M - f\opId + \opHess^h(f).
$$
is a trace-free Codazzi field, and existence follows. The second assertion follows by Stokes' Theorem (c.f. \cite{Tromba}), and this completes the proof.\qed
\medskip
Having established these preliminaries, we now prove Theorem \procref{MainTheorem}.
\medskip
{\noindent\bf Proof of Theorem \procref{MainTheorem}:}~The main work lies in proving surjectivity. Set $c_1:=\sqrt{-k_+}$ and $c_2:=\sqrt{-k_-}$. Let $h_1$ and $h_2$ be marked hyperbolic metrics in $\Thyp[S]$. Consider the function $\Psi:\Thyp[S]\rightarrow\Bbb{R}$ given by
$$
\Psi(h) := \Psi_{h_1,h_2}(h) := c_1F_{h_1}(h) + c_2F_{h_2}(h) = c_1F(h_1,h) + c_2F(h_2,h).\eqnum{\nexteqnno[DefinitionOfPhi]}
$$
By Lemma \procref{BonMonSch}, this function is proper and strictly convex with respect to the Weyl-Peterson metric. It therefore attains a unique minimum at some point $h_0$, say. Let $\rho$ be the point of $\Trep[S]$ corresponding to $h_0$. We will show that $\rho$ is the Fuchsian component of an affine deformation whose corresponding marked GHMC Minkowski spacetime $X$ solves
$$
\Phi(k_+,k_-,X)=((S,h_1),(S,h_2)).
$$
\par
By \eqnref{FormulaForDerivative}, at this point
$$
\int_\Sigma\opTr(A(c_1B(h_1,h_0)+c_2B(h_2,h_0)))\opdArea_{h_0} = 0,\eqnum{\nexteqnno[IntegrabilityCondition]}
$$
for every trace-free Codazzi field $A$ of $h_0$. It follows by Lemma \procref{OrthogonalDecomposition} that
$$
c_1B(h_1,h_0)+c_2B(h_2,h_0) = f\opId - \opHess_{h_0}(f),\eqnum{\nexteqnno[TheSumIsADerivative]}
$$
for some smooth function $f:S\rightarrow\Bbb{R}$.
\par
Now identify $S$ with the quotient of $\Bbb{H}^2$ by $\rho(\pi_1(S))$ and identify $B(h_1,h_0)$, $B(h_2,h_0)$ and $f$ with their lifts to $\Bbb{H}^2$. Next identify $\Bbb{H}^2$ with the future component of the unit pseudo-sphere in $\Bbb{R}^{2,1}$ and identify $T\Bbb{H}^2$ with a subbundle of the trivial bundle $\Bbb{H}^2\times\Bbb{R}^{2,1}$. In particular, this identifies the fields $B(h_1,h_0)$ and $B(h_2,h_0)$ with sections of $T^*\Bbb{H}^2\otimes\Bbb{R}^{2,1}$. Given a base point $x_0\in\Bbb{H}^2$, the functions $X_1,X_2:\Bbb{H}^2\rightarrow\Bbb{R}^{2,1}$ are defined by
$$\eqalign{
X_1 &:= U_1 + \int_{x_0}^xc_1B(h_1,h_0)\cdot\partial_\tau d\tau,\ \text{and}\cr
X_2 &:= -U_2-\int_{x_0}^xc_2B(h_2,h_0)\cdot\partial_\tau d\tau,\cr}\eqnum{\nexteqnno[DefinitionOfX]}
$$
where $U_1$ and $U_2$ are constant vectors to be determined. Indeed, since $B(h_1,h_0)$ and $B(h_2,h_0)$ are Codazzi tensors, both integrals are well-defined independent of the paths chosen from $x_0$ to $x$.
\par
The functions $X_1$ and $X_2$ define respectively future- and past-oriented strictly convex embeddings of constant extrinsic curvature equal to $k_+$ and $k_-$ respectively. Furthermore, these embeddings are equivariant with respect to the affine deformations $(\rho,\tau_1)$ and $(\rho,\tau_2)$ respectively, where the cocycles $\tau_1$ and $\tau_2$ are given by
$$\eqalign{
\tau_1(\gamma) &:=\int_{x_0}^{\gamma(x_0)}c_1B(h_1,h_0)\cdot\partial_\tau d\tau + U_1 - \rho(\gamma)U_1,\ \text{and}\cr
\tau_2(\gamma) &:=-\int_{x_0}^{\gamma(x_0)}c_2B(h_2,h_0)\cdot\partial_\tau d\tau - U_2 + \rho(\gamma)U_2.\cr}
$$
It remains to show that $\tau_1=\tau_2$. However, in the present context, the support functions of $X_1$ and $X_2$ are given respectively by
$$\eqalign{
\phi_1(x) &:= \langle X_1(x),x\rangle,\ \text{and}\cr
\phi_2(x) &:= -\langle X_2(x),x\rangle.\cr}
$$
and it follows upon integrating \eqnref{TheSumIsADerivative} that
$$
\phi_1(x) + \phi_2(x) = f(x) + \langle U_1 + U_2 - \nabla f(x_0),x\rangle.
$$
We therefore set
$$
U_1 + U_2 = \nabla f(x_0),
$$
so that, since $f$ is bounded,
$$
\mlim_{x\rightarrow\partial_\infty\Bbb{H}^2}\frac{\phi_1(x)}{\opCosh(d(x,x_0))} = -\mlim_{x\rightarrow\partial_\infty\Bbb{H}^2}\frac{\phi_2(x)}{\opCosh(d(x,x_0))},\eqnum{\nexteqnno[FunctionsCoincideAtInfinity]}
$$
where $d(x,x_0)$ here denotes the distance in $\Bbb{H}^2$ from $x_0$ to $x$. Recall now that $\phi_1$ and $\phi_2$ define continuous functions $\tilde{\phi}_1$ and $\tilde{\phi}_2$ over $\partial_\infty\Bbb{H}^2$ (c.f. \cite{BarbotFillastre}). Furthermore, the convex hulls in co-Minkowski space $\Bbb{H}^2\times\Bbb{R}$ of the graphs of these functions are the dual convex sets of the Mess' future and past domains of $(\rho,\tau_1)$ and $(\rho,\tau_2)$ respectively. By \eqnref{FunctionsCoincideAtInfinity}, $\tilde{\phi}_1=\tilde{\phi}_2$ so that these Mess' future and past domains are components of the same GHMC Minkowski spacetime. Since this uniquely defines the cocycle, we conclude that $\tau_1=\tau_2$, and surjectivity follows.
\par
Observe now that if $X$ is the GHMC Minkowski spacetime associated to the affine deformation $(\rho,\tau)$, if $h\in\Thyp[S]$ is the marked hyperbolic metric corresponding to $\rho$, and if $\Phi(k_+,k_-,X)=((S,h_1),(S,h_2))$, then $h$ is in fact the unique minimum of the function $\Psi_{h_1,h_2}$ defined above. This proves injectivity. Finally, by elliptic regularity, $\Phi$ is real-analytic, and since $\Psi$ is strictly convex, analyticity of $\Phi^{-1}$ follows by the inverse function theorem for analytic functions. This completes the proof.\qed
\global\subheadno=0
\inappendicestrue
\medskip
\newsubhead{Bibliography}[Bibliography]
\medskip
{\leftskip = 5ex \parindent = -5ex
\leavevmode\hbox to 4ex{\hfil \cite{AndEtAl}}\hskip 1ex{Andersson L., Barbot T., Benedetti R., Bonsante F., Goldman W. M., Labourie F., Scannell K. P., Schlenker J. M., Notes on: ``Lorentz spacetimes of constant curvature'', {\sl Geom. Dedicata}, {\bf 126}, (2007), 47--70}
\medskip
\leavevmode\hbox to 4ex{\hfil \cite{BarbotBeguinZeghib}}\hskip 1ex{Barbot T., B\'eguin F., Zeghib A., Prescribing Gauss curvature of surfaces in 3-dimensional spacetimes: application to the Minkowski problem in the Minkowski space, {\sl Ann. Inst. Fourier}, {\bf 61}, (2011), no. 2, 511--591}
\medskip
\leavevmode\hbox to 4ex{\hfil \cite{BarbotFillastre}}\hskip 1ex{Barbot T., Fillastre F., Quasi-Fuchsian co-Minkowski manifolds, arXiv:1801.10429}
\medskip
\leavevmode\hbox to 4ex{\hfil \cite{BonMonSch}}\hskip 1ex{Bonsante F., Mondello G., Schlenker J. M., A cyclic extension of the earthquake flow {II}, {\sl Ann. Sci. \'Ec. Norm. Sup\'er.}, {\bf 48}, (2015), no. 4, 811--859}
\medskip
\leavevmode\hbox to 4ex{\hfil \cite{Danciger}}\hskip 1ex{Danciger J., A Geometric transition from hyperbolic to anti de Sitter geometry, {\sl Geom. Topol.}, {\bf 17}, (2013), no. 5, 3077--3134}
\medskip
\leavevmode\hbox to 4ex{\hfil \cite{FillSeppi}}\hskip 1ex{Fillastre F., Seppi A., Spherical, hyperbolic and other projective geometries, in {\sl Eighteen essays on non-commutative geometry}, IRMA, {\sl Math. Theor. Phys.}, {\bf 29}, Eur. Math. Soc., Z\"urich, (2019)}
\medskip
\leavevmode\hbox to 4ex{\hfil \cite{FillSmi}}\hskip 1ex{Fillastre F., Smith G., Group actions and scattering problems in Teichm\"uller theory, in {\sl Handbook of group action IV}, Advanced Lectures in Mathematics, {\bf 40}, (2018), 359--417}
\medskip
\leavevmode\hbox to 4ex{\hfil \cite{Goldman}}\hskip 1ex{Goldman W., The symplectic nature of the fundamental group of surfaces, {\sl Adv. Math.}, {\bf 54}, (1984), no. 2, 200--225}
\medskip
\leavevmode\hbox to 4ex{\hfil \cite{LabourieIII}}\hskip 1ex{Labourie F., Probl\`eme de Minkowski et surfaces \`a courbure constante dans les vari\'et\'es hyperboliques, {\sl Bull. Soc. math. France}, {\bf 119}, 1991, 307--325}
\medskip
\leavevmode\hbox to 4ex{\hfil \cite{Labourie}}\hskip 1ex{Labourie F., Surfaces convexes dans l'espace hyperbolique et $\Bbb{CP}^1$-structures, {\sl J. London Math. Soc.}, {\bf 45}, (1992), no. 3, 549--565}
\medskip
\leavevmode\hbox to 4ex{\hfil \cite{LabourieII}}\hskip 1ex{Labourie F., M\'etriques prescrites sur le bord des vari\'et\'es hyperboliques de dimension $3$, {\sl J. Diff. Geom.}, {\bf 35}, (1992), no. 3, 609--626}
\medskip
\leavevmode\hbox to 4ex{\hfil \cite{Mess}}\hskip 1ex{Mess G., Lorentz spacetimes of constant curvature, {\sl Geom. Dedicata}, {\bf 126}, (2007), 3--45}
\medskip
\leavevmode\hbox to 4ex{\hfil \cite{Schlenker}}\hskip 1ex{Schlenker J. M., Hyperbolic manifolds with convex boundary, {\sl Inventiones math.}, {\bf 163}, (2006), 109--169}
\medskip
\leavevmode\hbox to 4ex{\hfil \cite{Tamburelli}}\hskip 1ex{Tamburelli A., Prescribing metrics on the boundary of anti-de Sitter $3$-manifolds, {\sl Int. Math. Res. Not. IMRN}, (2018), no. 5, 1281--1313}
\medskip
\leavevmode\hbox to 4ex{\hfil \cite{Tromba}}\hskip 1ex{Tromba A. J., {\sl Teichm\"uller theory in Riemannian geometry}, Lectures in Mathematics ETH Z\"urich, Birkh\"auser Verlag, Basel, (1992)}
\par}
%
%%%%%%%%%%%%%%%%%%%%%%%%%%%%%%%%%%%%%%%%%%%%%%%%%%%%%%%%%%%%%%%%%%%%%%%%%%%%%%%%%%%%%%%%%%%%%%%%%%%%%%%%%%%%%%%%%%%%%%%
%
% 3: Closing commands.
%
%%%%%%%%%%%%%%%%%%%%%%%%%%%%%%%%%%%%%%%%%%%%%%%%%%%%%%%%%%%%%%%%%%%%%%%%%%%%%%%%%%%%%%%%%%%%%%%%%%%%%%%%%%%%%%%%%%%%%%%
%
\enddocument